\documentclass[a4paper,11pt,ralign]{amsart}
\usepackage{amsmath,amssymb,mathabx,amsfonts,graphicx}
\usepackage[mathscr]{euscript}
\usepackage[all]{xy}
\usepackage{enumitem}
\usepackage{inputenc}
\usepackage{xcolor}
\setenumerate[1]{label={\textup{(\arabic*)}}}
    \theoremstyle{plain}
\newtheorem{theorem}{Theorem}[section]

\newtheorem{lemma}[theorem]{Lemma}
\newtheorem{proposition}[theorem]{Proposition}
\newtheorem{corollary}[theorem]{Corollary}
\newtheorem*{corollary*}{Corollary}

\newtheorem{definition}[theorem]{Definition}
\newtheorem{definitions}[theorem]{Definitions}

\theoremstyle{definition}
\newtheorem{remark}[theorem]{Remark}

\theoremstyle{remark}
\newtheorem{example}[theorem]{Example}

\newcommand{\G}{\ensuremath \mbox{\tiny $G$}}
\newcommand\m[1]{\mathrm{m}_{\mbox{\tiny \ensuremath  #1}}}	
\newcommand{\dmg}{\ensuremath  \mathrm{dm}_{\G}}
\newcommand{\mg}{\ensuremath  \mathrm{m}_{\G}}

\newcommand{\supp}{\ensuremath  \mathrm{supp}}

\newcommand\lE[1]{\ensuremath \mathcal{L}\left(#1\right)}
\newcommand{\mh}{\ensuremath  \, \mathrm{m}_{\mbox{\tiny $H_\mu$}}}
\newcommand\<[1]{\left\langle\, #1\,\right\rangle}

\newcommand\norm[1]{\lVert #1 \rVert}
\newcommand\bnorm[1]{\Big\lVert #1 \Big\rVert}

\newcommand{\Cf}{\ensuremath \mbox{\large${\mathbf{1}}$}}

\newcommand{\h}{\ensuremath \mathbb{H}}

\newcommand{\dmh}{\ensuremath  \, \mathrm{dm}_{\mbox{\tiny $H$}}}

\newcommand{\mA}{\mathscr{A}}

\newcommand\restr[2]{{
  \left.\kern-\nulldelimiterspace 
  #1 
  \vphantom{|} 
  \right|_{_{#2}} 
  }}

\newcommand{\C}{\mathbb{C}}
\newcommand{\D}{\mathbb{D}}
\newcommand{\T}{\mathbb{T}}

\newcommand{\N}{\mathbb{N}}

\newcommand{\Z}{\mathbb{Z}}

\newcommand{\sap}{\sigma_{ap}}
\newcommand{\sapT}{\sigma_{ap}(T)}
\newcommand{\sapuno}{\sigma_{ap}(\lambda_1(\mu))}
\newcommand{\sapunocero}{\sigma_{ap}(\lambda_1^0(\mu))}
\renewcommand{\emptyset}{\varnothing}

\DeclareMathOperator{\acc}{Acc}

\title[Uniformly ergodic probability measures]{Uniformly ergodic probability measures}

\author{Jorge Galindo,  Enrique Jord\'a \and Alberto Rodr\'iguez-Arenas}

\address{\noindent Jorge Galindo, Instituto Universitario de Matem\'aticas y
Aplicaciones (IMAC)\\ Universidad Jaume I, E-12071, Cas\-tell\'on,
Spain. \hfill\break \noindent E-mail: {\tt jgalindo@uji.es}}
\address{\noindent Enrique Jord\'a, EPS Alcoy, Instituto Universitario de Matem\'atica Pura y Aplicada IUMPA,
Universitat Polit\`ecnica de Val\`encia, Plaza Ferr\'andiz y Carbonell s/n
E-03801 Alcoy, Spain, \hfill\break \noindent E-mail: {\tt ejorda@mat.upv.es}}
\keywords{ergodic measure, uniformly ergodic measure, random walk, mean ergodic operator, uniformly mean ergodic operator, convolution operator, locally compact  group, measure algebra}
\address{\noindent Alberto Rodr\'{\i}guez-Arenas \\Departamento de Matem\'atiques, Universidad Jaume I, E-12071, Cas\-tell\'on,
Spain. \hfill\break \noindent E-mail: {\tt arenasa@uji.es} }
\keywords{ergodic measure, uniformly ergodic measure, random walk, mean ergodic operator, uniformly mean ergodic operator, convolution operator, locally compact  group, measure algebra}

\subjclass[2020]{Primary 43A05; Secondary 43A20, 43A30, 46H99, 47A35}

\date{\today.}

\begin{document}

\begin{abstract}
Let $G$ be a locally compact group and $\mu$ be a probability measure on $G$.  We consider the convolution operator $\lambda_1(\mu)\colon L_1(G)\to L_1(G)$ given by $\lambda_1(\mu)f=\mu \ast f$ and  its  restriction $\lambda_1^0(\mu)$  to the augmentation ideal $L_1^0(G)$. Say  that $\mu$ is uniformly ergodic if the Ces\`aro means of the  operator  $\lambda_1^0(\mu)$ converge uniformly to 0, that is, if   $\lambda_1^0(\mu)$ is a uniformly mean ergodic operator with limit 0, and  that $\mu$ is uniformly completely mixing if the powers of the  operator  $\lambda_1^0(\mu)$ converge uniformly to 0.

We completely characterize the uniform mean ergodicity, and the uniform convergence of its powers, of  the operator $\lambda_1(\mu)$  and see that there is no difference between $\lambda_1(\mu)$ and   $\lambda_1^0(\mu) $ in these regards. We prove in particular  that $\mu$ is uniformly ergodic if and only if $G$ is
compact, $\mu$ is adapted (its support is  not contained in a proper closed subgroup of $G$) and 1 is an isolated point of the spectrum of $\mu$. The last of these three   conditions can actually be replaced by  $\mu$ being spread-out (some convolution power of $\mu$ is not singular). The measure $\mu$ is uniformly completely mixing if and only if $G$ is compact,  $\mu$ is spread-out and the only unimodular value in the spectrum of $\mu$ is 1.

\end{abstract}

\maketitle

\section{Introduction}
A probability measure $\mu$  on  a locally compact group $G$ defines  the transition probabilities  of a random walk on $G$. It is common knowledge in Ergodic Theory that the ergodicity of  a random walk can be characterized    in terms of the  fixed points of its Markov operator.  In the present case, the Markov operator is  identified with the \emph{convolution operator} $\lambda_1(\mu)\colon L_1(G)\to L_1(G)$, $\lambda_1(\mu)f=\mu \ast f$. This  is reflected  in the following theorem  stated in Rosenblatt \cite{rose81}  that can be regarded as the starting point of this work.  \begin{theorem}[Proposition 1.2 of \cite{rose81}]\label{rose81}
   Let $G$ be a locally compact group and let $\mu $ be a probability measure on $G$. Consider the operator $\lambda_1^0(\mu) $ that results from restricting $\lambda_1(\mu)$ to the augmentation ideal  \[L_1^0(G)=\left\{f\in L_1(G)\colon \int f(x)\dmg(x)=0\right\}.\] The following are equivalent:
    \begin{enumerate}
      \item The random walk defined by  $\mu$ is ergodic.
      \item The Ces\`aro means of the operator  $\lambda_1^0(\mu) $ converge to 0 in the strong operator topology.
          \item If $f\in L_\infty(G)$ and $\mu^\ast \ast f=f$, where $\mu^\ast (A)=\mu(A^{-1})$, then $f$ is constant almost everywhere.
    \end{enumerate}
 \end{theorem}
A probability measure  $\mu$ can be ergodic (and even spread-out) while $\mu^\ast $ is not;  see
Remark on p. 9 of S. Glasner [12]   which leans on  results of Azencott \cite{azen70}. This shows that we cannot use $\mu$ instead of $\mu^*$ in condition (3) of Theorem \ref{rose81}.

If the group $G$ is abelian or compact, the Choquet-Deny  and Kawada-It\^o theorems, \cite{choqdeny60,kawaito}, show, respectively, that   the random walk induced by $\mu$  is ergodic (and we then say that $\mu$ itself is \emph{ergodic}) if and only if its support is not contained in a proper closed subgroup of $G$, i.e, if $\mu$ is \emph{adapted}.
If  the support of $\mu$   is not even contained in a translate of a proper closed normal subgroup (we say then that $\mu$ is \emph{strictly aperiodic}),  $\mu$ has  a stronger property, the strong operator limit of $\lambda_1^0(\mu^n)$ is 0. We say then that  $\mu$ is \emph{completely mixing}. All these facts are well-known and collected in Theorem \ref{adp-spbasico}, below.

In this paper we address the problem of characterizing ergodicity when
 the strong operator topology is replaced by the operator norm. Our approach will be operator-theoretic and focuses on  the uniform mean ergodicity (operator norm convergence of the Ces\`aro means) of the operators $\lambda_1(\mu)$ and $\lambda_1^0(\mu)$   and on the uniform ergodicity of $\mu$, where $\mu$ is said to be  \emph{uniformly ergodic} if   ${\displaystyle \lim_n \bnorm{\dfrac{1}{n}\sum_{k=1}^n\lambda_1^0(\mu)^k}=0}$.

After  two preliminary sections, we devote Section 3 to  the operator $\lambda_1^0(\mu)$ and its  relation with  uniform ergodicity and
 Section 4 to present a first  characterization of uniform ergodicity which shows its several angles and  reflects what was previously known on the subject.

     In Section 5  spectral conditions are introduced. It is shown in \cite{galijorda21} that, for an adapted measure $\mu$ and an abelian group $G$,  $\lambda_1(\mu)$ is uniformly mean ergodic if and only if the support of $\mu$ is contained in a compact subgroup of $G$  and $1$ is isolated  in the spectrum of $\mu$. One of our main results here, Theorem \ref{car:adap}, is that this theorem still holds true when $G$ is not abelian, and that the condition of having $1$ as an  isolated point of the spectrum of $\mu$  is in fact equivalent,  when $G$ is compact and $\mu$ is adapted,   to $\mu$ being spread-out.

The subject of uniform mean ergodicity of convolution operators has been discussed in the book by  Revuz, \cite[Section 6.3]{revuz}, where uniformly ergodic convolution operators on $L_\infty(G)$ that converge to a rank-one operator are characterized, as well as in a brief incursion by Mustafayev and Topal \cite[Theorem 1.2]{mustatopal21}. To the best of our knowledge, no spectral characterizations of uniform ergodicity of convolution operators  on $L_1(G)$, $G$ any locally compact group, have been obtained before.

We finally complete our study of the operator $\lambda_1^0(\mu)$ and characterize quasi-compactness and convergence of the powers  of the operators $\lambda_1(\mu)$ and $\lambda^0_1(\mu)$. Convergence in norm of the sequence of iterates $(\lambda_1^0(\mu^n))$ is shown to be equivalent to that of $(\lambda_1(\mu^n))$,  providing   a spectral extension of results of   Anoussis and Gatzouras in \cite{anougatz04}.
\section{Preliminaries}\label{sec:prel}
We use this section to set some  notations and recall  a few  definitions and basic known facts.
\subsection{Mean ergodic operators}
For a bounded linear operator $T:X\longrightarrow X$ on a Banach space $X$, we use the following notations for the iterates and averages:
\[
T^n = T\circ \overset{n}{\cdots} \circ T, \quad T_{_{[n]}} = \frac{T+\cdots + T^n}{n}.
\]
\begin{definitions}
Let $X$ be a Banach space and let $T:X\longrightarrow X$ be a bounded linear operator. The operator is
\begin{itemize}
    \item [(i)] \emph{mean ergodic} if the sequence  $(T_{_{[n]}}x)_n $ is convergent for every $x\in X$.
    \item [(ii)] \emph{uniformly mean ergodic}  if the sequence  $(T_{_{[n]}})_n $ is       convergent in the operator norm.
\end{itemize}
\end{definitions}
If $T$ is mean ergodic, the limit of $(T_{_{[n]}}x)_n $  is  always  a linear  projection $P:X\longrightarrow X$, satisfying $\ker P = \overline{(I-T)X}$. Moreover,  $X=\ker P \oplus \ker (I-T)$.

Uniform mean ergodicity for power bounded operators is completely characterized by Lin.
\begin{theorem}[Theorem of \cite{lin74}]
\label{lin}
Let $T$ be a power bounded (i.e., so that $(\norm{T^n})_n$ is a bounded sequence) linear operator on a Banach space $X$.  The operator $T$ is uniformly mean ergodic if, and only if,  $(I-T)X$ is closed. The Banach space $X$ can be then decomposed as $X=(I-T)X\bigoplus \ker(I-T)$.
\end{theorem}

Lin's theorem completed the  following theorem, due to  Dunford \cite{dunf43},  which reveals that the  uniform mean ergodicity  of an operator  can be characterized through its spectral properties. See also \cite[Theorem 2.2.7]{kren85}
 and Lin \cite{lin74}.

\begin{theorem}
\label{dunfordlin}
Let $T$ be a power bounded linear operator on a complex Banach space $X$. Then $T$ is uniformly mean ergodic if, and only if, either $1\not\in \sigma(T)$ or it is a pole of order $1$ of the resolvent.
\end{theorem}

For the next result, due to Yosida and Kakutani, we recall that an operator $T$ is quasi-compact if there exists a compact operator $K$ such that $\|T^n-K\|<1$ for some $n\geq 1$. By $\sigma_p(T)$, we will  denote the point-spectrum of $T$, the set of its eigenvalues, and $\T$ will denote the set of complex numbers of modulus 1.
\begin{theorem}[Theorem 4 and its Corollary in page 205 of \cite{yosikaku41}]
\label{yosidakaku}
 Let $T$ be a power bounded and quasi-compact linear operator on a complex Banach space $X$. The following assertions hold:
 \begin{itemize}
     \item[(i)] $(T_{_{[n]}})_n$ converges uniformly to a finite rank projection $P$.
     \item[(ii)] $(T^n)_n$ converges in norm to a finite rank projection $P$ if, and only if,  $\sigma_p (T)\cap \T \subseteq \{1\}$.
     \item[(iii)] $(T^n)_n$ converges in norm to $0$ if, and only if, $\sigma_p (T)\cap \T =\emptyset$.
 \end{itemize}
 \end{theorem}

\subsection{Convolution operators}
We will be working with convolutions  on $\sigma$-compact locally compact groups. If $G$ is a locally compact, $\sigma$-compact group, $e$ will stand for its identity element and  $\mg$  for its Haar measure, the (essentially unique) Borel measure that is left invariant on a locally compact group $G$. The Banach space of continuous functions on $G$ vanishing at infinity is denoted by $C_0(G)$ and that of compactly supported functions, by $C_{00}(G)$. The Banach space of bounded regular measures is $M(G)=C_0(G)^*$. We denote by $L_1(G,\mu)$ the Banach space of equivalence classes of integrable functions with respect to a measure $\mu$ on $G$. If $\mu=\mg$ is the Haar measure, we simply write $L_1(G)$. We denote by $L_1^0(G)$ the augmentation ideal, the space of functions with integral $0$. We regard $L_1(G)$ as an ideal of $M(G)$ through the embedding $f\mapsto f\cdot d\mg$.

\begin{definitions}
Let $\mu,\mu_1,\mu_2 \in M(G)$. We consider
\begin{itemize}
    \item[(i)] the convolution of measures:
    \[
    \langle\mu_1*\mu_2,h \rangle=\int \int h(xy) d\mu_1 (x) d\mu_2 (y), \quad \text{for $h \in C_{00}(G)$}.
    \]
    \item[(ii)] the left convolution operators: $\lambda_p(\mu):L_p(G)\longrightarrow L_p(G)$, $1\leq p\leq \infty$, given by
    \[
        \lambda_p (\mu)f(s) =(\mu*f)(s)=\int f(x^{-1}s) d\mu (x), \quad \text{for } f\in L_p(G), s\in G.
    \]
    \item[(iii)] the  operator $\lambda_1 ^0 (\mu)\colon L_1^0(G)\to L_1^0(G)$ that arises when $\lambda_1(\mu)$ is restricted to the augmentation ideal $L_1 ^0 (G)$.
    \item[(iv)] the convolution iterates  $\mu ^n$ and the convolution averages $\mu_{_{[n]}}$, respectively: $
        \mu^n = \mu *\cdots *\mu$ and $ \mu_{_{[n]}} = \frac{\mu + \cdots +\mu^n}{n}.$
\end{itemize}
\end{definitions}
\subsection{Ergodic measures}
A measure $\mu \in M(G)$ on a locally compact group $G$ is said to be ergodic when the random walk it generates is ergodic and it is said to be completely mixing when the random walk is mixing. For an operator-theoretic approach, however,  the following, equivalent, definitions are more suitable. They  are, in fact, the ones chosen by most authors, see, e.g., \cite{anougatz08,cuny03,jawo04,linwitt94,rose81} (note that Rosenblatt \cite{rose81} uses the terms \emph{ergodic by convolutions} and \emph{mixing by convolutions}).

\begin{definitions}\label{def:measerg+cm}
A probability measure $\mu\in M(G)$ is:
\begin{itemize}
    \item[\normalfont (i)] \emph{ergodic} if $\lim_n \|\mu_{_{[n]}}*f\|_1 = 0 $, for all $f\in L_1^0(G)$. In other words, if $\lambda_1 ^0 (\mu)$ is  mean ergodic and the strong operator limit of the means is the null operator.
    \item[\normalfont (ii)] \emph{completely mixing} if the sequence of iterates $(\lambda_1^0(\mu^n)) _{n\in\N}$ converges to 0 in the strong operator topology.
   \end{itemize}
\end{definitions}
%

Azencott \cite[p. 43, Corollaire]{azen70}
proved that the existence of an ergodic probability implies amenability (see also \cite[Proposition 1.9]{rose81}).  J. Rosenblatt \cite[Theorem 1.10]{rose81}
 proved that in any amenable group there is an ergodic probability.

 If  in  Definition \ref{def:measerg+cm} the strong operator topology is replaced by the operator norm, two new concepts arise.
\begin{definitions}
  We say that  probability measure $\mu\in M(G)$ is:
\begin{itemize}
    \item[(i)]  \emph{uniformly ergodic} if $\lambda_1 ^0 (\mu)$ is uniformly mean ergodic and the operator norm limit of the means  is the null operator.
     \item[(ii)] \emph{uniformly completely mixing} if the sequence of iterates $(\lambda_1^0(\mu^n)) _{n\in\N}$ converges to 0 in the operator norm.
\end{itemize}
\end{definitions}

Note that uniform ergodicity of $\mu\in M(G)$  comprises convergence to 0 of the means of the operator $\lambda_1^0(\mu)$.  The example  $\mu =\delta_e$ shows that this property is  stronger  than uniform mean ergodicity of $\lambda_1^0(\mu)$.

The analysis of the (uniform) ergodicity or the (uniform) complete mixing of a measure $\mu$ is heavily related to the algebraic properties of its support subgroup $H_\mu$.  We denote by $H_\mu$ the smallest closed subgroup of $G$ that contains $S_\mu$, the support of $\mu$.
 To give a hint on this relation we need to set some  definitions that  have already been informally introduced.
\begin{definitions}
A probability measure $\mu\in M(G)$ is:
\begin{itemize}
    \item[(i)] \emph{adapted} if $H_\mu =G$.
    \item[(ii)] \emph{strictly aperiodic} if it is adapted and the only normal closed subgroup $N$ satisfying $S_\mu \subseteq xN$ for some $x\in G$ is actually $N=G$.
    \item[(iii)] \emph{spread-out} if $\mu ^n $ is not singular for some $n\in \N$, i.e. if there is $n\in \N$ such that  $\mu^n (A) =1$ implies $\mg (A)>0$.
\end{itemize}
\end{definitions}
We collect now some  well-known facts that will be needed later on. The first two statements  are proved by Lin and Wittman in Corollary 2.7 of \cite{linwitt94}. The third statement  is the  Kawada-It\^o theorem \cite{kawaito}. For the last sentence one has to apply   Proposition 3.1 of \cite{galijorda21} to the previous one. See also Corollary 3.2 of \cite{derrlin89}.
\begin{theorem} \label{adp-spbasico}
Let $G$ be a locally compact group and let $\mu\in M(G)$ be a probability measure. The following assertions hold
\begin{enumerate}
    \item If $\mu$ is ergodic, then  $G=\overline{\bigcup_j \bigcup_k S_\mu ^{-j} S_\mu ^k}$. Hence $\mu$ ergodic implies that $\mu$ is adapted.
        \item    If $\mu$ is completely mixing, then $G=\overline{\bigcup_j S_\mu ^{-j} S_\mu ^j}$. Hence, $\mu$  completely mixing implies that  $\mu$ is strictly aperiodic.
     \item If $G$ is compact, $(\mu_{[n]})_{n\in\N}$ converges, in the weak* topology of $M(G)$, to the Haar measure of $H_\mu$.  If $\mu$ is strictly aperiodic,  $(\mu^n)_{n\in\N}$ converges, in the weak* topology of $M(G)$, to the Haar measure of $G$. The measure  $\mu$ is then  completely mixing.
     \end{enumerate}
\end{theorem}

\begin{remark}
Complete mixing and ergodicity  are   known to be equivalent for strictly aperiodic measures under some restrictions on the measure and/or the group.  This is the case, for instance  when the group is in [SIN] (i.e. has a neighbourhood basis consisting of  sets $U$ such that $Ux=xU$ for every $x\in G$), see \cite[Theorem 5.2]{jawo04}. Looking from  the measure side,  Glasner \cite{glas76} (leaning on Foguel \cite{fogu75b})   showed that complete mixing and ergodicity are equivalent for any spread-out strictly aperiodic measure.   Some other cases have been sorted out in \cite{anougatz08} and  \cite{cuny03} but   the problem of the equivalence, \emph{the complete mixing problem}, remains open.
\end{remark}
\subsection{Fourier-Stieltjes transforms}
To study the spectrum of $\mu$ when $G$ is compact, we will rely  on the Fourier-Stieltjes transform which associates   every continuous irreducible unitary representation $\pi$  of $G$ on a Hilbert space $\h_\pi$  with a homomorphism between $M(G)$ and $\mathcal{L}(\h_{\pi})$ the bounded operators on $\h_\pi$.
By a  unitary representation of $G$, we understand a continuous homomorphism of $G$ into the  group of unitary operators on a finite dimensional Hilbert space $\h_\pi$, where the latter is assumed to carry the weak operator topology.
 A representation $\pi$ of $G$ on a Hilbert space $\h$ is \emph{irreducible} if there is no nontrivial closed subspace $V$ of $\h$ such that $\pi(G)V\subseteq V$.
\begin{definition}[Fourier-Stieltjes transform]
  If  $\mu\in M(G)$ is a bounded measure on a compact group $G$ and $\pi\colon G \to \mathcal{U}(\mathbb{H}_\pi)$ is a unitary representation of $G$ on a finite dimensional Hilbert space $\mathbb{H}_\pi$, the    \emph{Fourier-Stieltjes transform} of $\mu$ at $\pi$ is the operator $\widehat{\mu}(\pi)\in \lE{\mathbb{H}_\pi}$ defined by
  \[ \<{\widehat{\mu}(\pi) \xi, \eta}=\int_{\G} \overline{\<{\pi(t)\xi,\eta}} d\mu(t),\quad \mbox{ for every $\xi,\eta\in \mathbb{H}_\pi$}.\]
\end{definition}

The basic properties of Fourier-Stieltjes transforms are collected  below. Statements \ref{en:1}, \ref{en:2}, \ref{en:4} and \ref{en:6}  are proved in  Section 28 of \cite{hewiross2}. Statement \ref{en:5} is  Theorem 27.17 of \cite{hewiross2}.
 \begin{theorem}[Some properties of the Fourier-Stieltjes transform]\label{FS}
 Let $G$ be a compact group.
 \begin{enumerate}
   \item \label{en:1} If $\mu_i\in M(G)$, $i=1,2$, are bounded measures and $\pi\colon G \to \mathcal{U}(\mathbb{H}_\pi)$ is a unitary representation of $G$ on a finite dimensional Hilbert space $\mathbb{H}_\pi$, then $\widehat{\mu_1\ast \mu_2}(\pi)=\widehat{\mu_1}(\pi) \widehat{ \mu_2}(\pi)$.
   \item \label{en:2} If $ \mu \neq 0$, there is an irreducible unitary representation of $G$ on a finite dimensional Hilbert space, $\pi\colon G \to \mathcal{U}(\mathbb{H}_\pi)$ such that
       \[ \widehat{\mu}(\pi)\neq 0.\]
      \item \label{en:4} If $\mu\in M(G)$ and $\pi$ is an irreducible unitary representation of $G$ on a finite dimensional Hilbert space, then
       \[\sigma\left(\widehat{\mu}(\pi) \right)\subseteq \sigma(\mu).\]
   \item \label{en:6} Let   $\mathbf{1}$ denote the trivial 1-dimensional representation. If $\mathbf{1}\neq \pi$ is an irreducible representation of $G$, then $\widehat{\mg}(\pi)=0$.
     \item \label{en:5}  Let $\pi\colon G\to \h_\pi$ be an irreducible unitary representation. If $\{\xi_1,\ldots,\xi_n\}$  is  an orthonormal basis of $\h_\pi$, then
    \begin{align*}
     \int\left| \<{\pi(t)\xi_i,\xi_j}\right|^2 \dmg(t)&=\begin{cases} \frac{1}{n},\; \mbox{ if $i=j$, }\\ 0, \; \mbox{ if $i\neq j$ }.\end{cases}
   \end{align*}
        \end{enumerate}\end{theorem}
We add next two simple known facts that  will be  needed later on. We provide proofs for the sake of completeness.
\begin{corollary}
\label{en:3}
Let $G$ be a compact group and $\mu\in M(G) $ a bounded measure.
\begin{itemize}
\item[(a)]  For every irreducible representation $\pi$ and every $\xi\in \mathbb{H}_\pi$, $\|\xi\|=1$, there is $f_{\pi,\xi}\in L_1(G)$ such that  $ \widehat{f_{\pi,\xi}}(\pi)\xi=\xi$.
\item[(b)] If $\mu$  is adapted and $\pi$ is an irreducible unitary  representation of $G$ with $1\in \sigma\left(\widehat{\mu}(\pi)\right)$, then $\pi=\mathbf{1}$, the one-dimensional trivial representation.
 \end{itemize}
\end{corollary}
 \begin{proof}
We first  prove  (a).
  Define $f_{\pi,\xi} $ by
\[f_{\pi,\xi} (t)=n\<{\pi(t)\xi,\xi},\]
where $n=\dim \h_\pi$.
  Since $\pi$ is continuous for the weak operator topology, $f_{\pi,\xi}\in L_1(G)$.
  If $\{\xi,\xi_2,\ldots,\xi_n\}$ is an orthonormal basis of $\h_\pi$ containing $\xi$,
 Statement \ref{en:5} of Theorem \ref{FS}, shows that
  $ \<{\widehat{f_{\pi,\xi}}(\pi)\xi,\xi}=1$ and that  $ \<{\widehat{f_{\pi,\xi}}(\pi)\xi,\xi_j}=0$, for every $j=2,\ldots,n$. It follows that  $\widehat{f_{\pi,\xi}}(\pi)\xi=\xi$.

The proof of (b) is standard. Assume that $\mu$ is an  adapted measure and  that $\pi$ is  an irreducible representation  of $G$ with  $1\in  \sigma\left(\widehat{\mu}(\pi)\right)$. There is therefore $\xi \in \h_\pi$, $\norm{\xi}=1$, with
  \[ 1=\int\overline{\<{\pi(t)\xi,\xi}}\, d\mu(t).\]
  Since $|\<{\pi(t)\xi,\xi}|\leq 1$, for every $t\in G$, it follows that $ \<{\pi(t)\xi,\xi}=1$  for every $t\in S_\mu$ , and so that $\pi(t)\xi=\xi$ for every $t\in S_\mu$.
  We have then that \[S_\mu\subseteq  N_\xi:=\{t\in G\colon \pi(t)\xi=\xi\}.\]
  Since $N_\xi$ is a closed subgroup of $G$, adaptedness of $\mu$ implies that $N_\xi=G$.  As $\h_\pi$  cannot contain proper subspaces invariant under $\pi$, we conclude that $\pi=\mathbf{1}$.
\end{proof}\section{The norm of  $\lambda_1^0 (\mu)$}
The ergodicity properties of the measure $\mu$  derive  from those of the operator $\lambda_1^0(\mu)$.
While the norm of $\lambda_1(\mu)$ is always $\norm{\mu}$, see \cite[Page 47]{green69},   that is not necessarily the case for   $\lambda_1^0(\mu)$. It suffices to observe that, when $G$ is a compact group, $\lambda_1^0(\mg)=0$ . We see in this section that invariance of $\mu$ (and, \emph{a fortiori}, compactness of $G$) is the only obstacle against the equality $\norm{\lambda_1^0(\mu)}=\norm{\mu}$.
\begin{lemma}
\label{prop:apro}
Let $G$ be a locally compact group and let $\mu \in M(G)$ satisfy $\|\mu\|\leq 1$.  Then
      \[\norm{\lambda_1^0(\mu)}\geq  \frac{1}{ 2}\sup_{x\in G}\norm{\mu- \mu\ast \delta_x }.\]
\end{lemma}
\begin{proof}
Consider a two-sided approximate identity  $\{e_\alpha\colon \alpha\in \Lambda\}$ in $L_1(G)$,  with $\norm{e_\alpha}=1$ for every $
 \alpha\in \Lambda$, see \cite[Proposition 9.1.8]{dalesetal}. It is then easy to prove that $\{e_\alpha\colon \alpha\in \Lambda\}$ converges to $\delta_e$ in the $\sigma(M(G),C_0(G))$-topology. Since  convolution is separately $\sigma(M(G),C_0(G))$-continuous, we have that $(\mu- \mu\ast \delta_x)\ast e_\alpha$ converges to $\mu - \mu\ast \delta_x$ in the $\sigma(M(G),C_0(G))$-topology.

  This Proposition is actually based on the simple equality,
\begin{equation}\label{eq:apro}
    \mu\ast(e_\alpha-\delta_x \ast e_\alpha)=(\mu-\mu\ast \delta_x)\ast e_\alpha \quad \hfill   \mbox{($x\in G$  and $ \alpha \in \Lambda$)}.
\end{equation}

One can first use functions  $f_n\in C_0(G)$ so that
$ \norm{\mu-\mu\ast\delta_x}\leq \left|\<{\mu-\mu\ast\delta_x,f_n}\right|+1/(2n)
$  to show that,  for each $n\in \N$, there is $\alpha_n \in \Lambda$ with
\begin{equation*}
    \norm{\mu -\mu \ast \delta_x}\leq \norm{(\mu- \mu\ast \delta_x)\ast e_{\alpha_n}}+\frac{1}{n}.
\end{equation*}
From this, using \eqref{eq:apro} and taking into account that $e_\alpha-\delta_x\ast e_\alpha\in L_1^0(G)$,  we have
\begin{align*}
    \norm{\mu -\mu \ast \delta_x} &
\leq \norm{\lambda _1 ^0 (\mu)}\norm{e_{\alpha_n} - \delta_x \ast e_{\alpha_n}}+\frac{1}{n}, \text{ for each } n\in\N
\end{align*}
and hence
\[ \norm{\mu -\mu \ast \delta_x} \leq  2\norm{\lambda_1^0(\mu)},\]
yielding the result.
\end{proof}

\begin{proposition}
 \label{noncnorm}
      Let $G$ be a locally compact group that is \emph{not compact} and let $\mu\in M(G)$ with $\|\mu\|= 1$. Then
      $\norm{\lambda_1^0(\mu)}=1$.
\end{proposition}
    \begin{proof}
    We first assume that $S_\mu$ is  compact. Then there is $x\in G$ such that $S_\mu \cap S_\mu x=\emptyset$ which implies  that $\|\mu-\delta_x \ast \mu\|=2$.  We can apply then Lemma~\ref{prop:apro} and obtain that $\norm{\lambda_1^0(\mu)}=1$.

    In general, we can always find a sequence $\mu_n$ of measures with compact support such that $\mu_n$ converges in norm to $\mu$.
     Since $\norm{\lambda_1^0(\mu_n)}=1$ for every $n$,  $\norm{\lambda_1^0(\mu)}$ must be 1 as well.
    \end{proof}

\begin{corollary}
  \label{cor:noncnorm}
  Let $G$ be a locally compact group and let $\mu$ be a probability measure. If $\mu$ is uniformly ergodic, then $G$ is compact.
\end{corollary}
\begin{proof}
  The measure $\mu_{[n]}$ is again a probability measure. If $G$ is not compact, Proposition \ref{noncnorm} shows that $\norm{\lambda_1^0(\mu_{[n]})}=1$ and,  hence, $\mu$ cannot be uniformly ergodic.
\end{proof}
The previous Corollary could have  also been deduced from the following Proposition which should be compared to  Theorem \ref{adp-spbasico} (1).
\begin{proposition}
\label{caractporconjuntos}
  Let $\mu\in M(G)$ be a probability measure with support $S_\mu$.
 If $\mu$ is uniformly ergodic, then there exists $n\in \N$ such that \[G=\bigcup_{1\leq j, k\leq n} S_\mu ^{-j} S_\mu ^k.\]

\end{proposition}
\begin{proof}
 We proceed by contradiction. Assume $G\neq \bigcup_{1\leq j, k\leq n} S_\mu ^{-j} S_\mu ^k$, for any $n\in \N$.

  Since $S_{\mu_{_{[n]}}} = \bigcup_{j=1} ^n S_\mu ^j$, the assumption implies that $G\neq S_{\mu_{_{[n]}}}^{-1}S_{\mu_{_{[n]}}}$ and therefore, that,  for each $n\in\N$, there exists $x_n \in G$ such that $ S_{\mu_{_{[n]}}}\cdot x_n \cap S_{\mu_{_{[n]}}}=\emptyset$. From this, we have that
        $\|\mu_{_{[n]}}- \mu_{_{[n]}}\ast\delta_{x_n}\|=2$, for every $n\in\N$, and we conclude, upon applying Lemma~\ref{prop:apro},  that $\mu$ is not uniformly ergodic.
\end{proof}

\section{Uniform mean ergodicity of convolution operators}\label{sec:firstchar}
We approach in this section a characterization (Theorem \ref{firstcar:adap}) of uniform mean ergodicity of the operators $\lambda_1(\mu)$ and $\lambda_{1}^0(\mu)$ in terms of how close $\mu$ is to $\mh$.

Our attention is attracted right away to measures with $H_\mu$ compact. This is because of the following result proved in \cite[Theorem 5.4]{galijorda21}.

\begin{theorem}
\label{JOTL1}
Let $G$ be a locally compact group and let $\mu\in M(G)$ be a probability. Then $\lambda_1(\mu)$ is mean ergodic if and only if $H_\mu$ is compact.
\end{theorem}

We  first see that we can restrict our work to adapted measures (and compact groups).
\subsection{Reduction to the adapted case}
Let $G$ be  a locally compact group and let $ \mu\in M(G)$ be a bounded measure. We will use the symbol $\underline{\mu}$ to denote the measure $\mu$ when  seen as a measure in $M(H_\mu)$. We see here that $\lambda_1(\mu)$ is uniformly mean ergodic if and only if $\lambda_1(\underline{\mu})$ is.

We first note the relation between  uniform  convergence of powers and Ces\`aro means of the convolution operator $\lambda_1(\mu)$ with the weak$^\ast$-convergence of convolution powers and  Ces\`aro means of the measure itself.
  \begin{proposition}
  \label{ume0}
    Let $G$ be a locally compact group and let $\mu \in M(G)$ be a probability measure.
    \begin{itemize}
     \item[(a)] The  operator $\lambda_1(\mu)$ is uniformly mean ergodic if and only if $H_\mu$ is compact and
    \[ \lim_{n}\norm{ \lambda_1(\mu_{[n]}-\mathrm{m}_{\mbox{\tiny ${H_\mu}$}})}=0.\]
     \item[(b)] The sequence $(\lambda_1(\mu^n))$ is norm convergent if and only if $H_\mu$ is compact and
       \[ \lim_{n}\norm{ \lambda_1(\mu^n-\mathrm{m}_{\mbox{\tiny ${H_\mu}$}})}=0\]
     \end{itemize}
    \end{proposition}
    \begin{proof}
By Theorem \ref{JOTL1}, $H_\mu$ must be compact whenever $\lambda_1(\mu)$ is uniformly mean ergodic.

Theorem \ref{adp-spbasico} (3) shows that  $\mu_{[n]}$ converges to $\mathrm{m}_{\mbox{\tiny $H_\mu$}}$ in the $\sigma(M(H_\mu),C_0(H_\mu))$-topology. This  implies right away convergence in the $\sigma(M(G),C_0(G))$-topology.
The  weak$^\ast$-SOT sequential continuity of  $\lambda_1$, see \cite[Proposition 3.1]{galijorda21}, then shows that
 once we assume that either $\lambda_1(\mu_{[n]})$ or  $\lambda_1(\mu^n)$ is convergent, the limit must be to $\lambda_1(\mh)$.
    \end{proof}
    The desired equivalence is now  an immediate consequence of the last Proposition,  after recalling that $   \norm{\lambda_1(\nu)}=\norm{\nu}$, for every $\nu\in M(G)$, see \cite[Page 47]{green69}.
\begin{corollary}
  \label{3.4}
Let $G$ be a locally compact group and let $\mu\in M(G)$. The operator $\lambda_1(\mu)$ is uniformly mean ergodic if and only if $\lambda_1(\underline{\mu})$ is. Furthermore, $(\lambda_1(\mu^n))$ is norm convergent in $\lE{L_1(G)}$ if and only if $(\lambda_1(\underline{\mu}^n))$ is norm convergent in $\lE{L_1(H_\mu)}$.
 \end{corollary}
\subsection{Uniform mean ergodicity of $\lambda_1(\mu)$ vs uniform mean ergodicity of $\lambda_1^0(\mu)$}
Theorem \ref{JOTL1} shows that the operators  $\lambda_1^0(\mu)$ and   $\lambda_1(\mu)$ are very different from the point of view of mean ergodicity. If $G$ is abelian and $\mu\in M(G)$ is adapted, the operator $\lambda_1^0(\mu)$ is always mean ergodic (this is the Choquet-Deny theorem, \cite{choqdeny60}) but $\lambda_1(\mu) $ will only be mean ergodic when $G$  is compact, the case  $\mu=\delta_1\in M(\Z)$   is a particularly simple case.
The situation is completely different   in the uniform case.

\begin{proposition}
  \label{Rest} Let $T\in \lE{X}$ be a bounded operator on a Banach  space $X$ and let $X_0$ be a hyperplane of $X$ which is invariant under $T$, i.e., with $T(X_0)\subseteq  X_0$. Then $T$ is uniformly mean ergodic if and only $\restr{T}{X_0}$ is.
\end{proposition}
\begin{proof}
  It is clear that $T$ uniformly mean ergodic implies that $\restr{T}{X_0}$ is uniformly mean ergodic.

  Assume now that $\restr{T}{X_0}$ is uniformly mean ergodic. The ergodic decomposition of Theorem \ref{lin} yields then
  \[ X_0= (I-T)(X_0)\bigoplus \ker\left(I-\restr{T}{X_0}\right).\]
  Since $X_0$ is a hyperplane, there is $y\in X$ such that $X=X_0\oplus \<{y}$.
  This shows that $X=(I-T)(X_0)\oplus \ker\left(I-\restr{T}{X_0}\right)\oplus \<{y}$.

  Let now $a, b\in X_0$ and $\alpha \in \C$, be  such that $b=Tb$ and
  \[(I-T)y=(I-T)a+b+\alpha y,\]
  and define $\tilde{y}=b+\alpha y$. Then $\tilde{y}\in (I-T)(X)$. We now check that  $(I-T)(X)=(I-T)(X_0)\oplus \<{\tilde{y}}$.
  Pick to that end an arbitrary $x=x_0+\beta y\in X$, $x_0\in X_0$, $\beta \in \C$. Then
  \begin{align*}
    (I-T)x&=(I-T)x_0+\beta (I-T) y \\
    &=(I-T)x_0 + \beta (I-T)a + \beta (b+ \alpha y)  \\
   & \in (I-T)X_0 +\<{\tilde{y}}.
  \end{align*}
  Having proved that $(I-T)(X)=(I-T)(X_0)\oplus \<{\tilde{y}}$, we see that $(I-T)(X)$ is closed and, hence, that $T$ is uniformly mean ergodic by Theorem \ref{lin}.
\end{proof}
After noting that $L_1^0 (G)$, being the kernel of the integral functional, is a hyperplane of $L_1(G)$ that is invariant under $\lambda_1(\mu)$, the following Corollary is immediately obtained.
\begin{corollary} \label{cor:Rest}
  Let $G$ be a locally compact group and $\mu\in M(G)$ be a probability measure.  The operator $\lambda_1(\mu)$ is uniformly mean ergodic if and only if $\lambda_1^0(\mu)$ is uniformly mean ergodic.
\end{corollary}
\begin{remark} $\,$\ \begin{itemize} \item[(a)] The proof of Proposition \ref{Rest} shows that  the equivalence between the uniform mean ergodicity of $\lambda_1(\mu)$ and that of $\lambda_1^0(\mu)$ holds for any measure $\mu \in M(G)$, with $\|\mu^n\|\leq M,$ $n\in \N$.
	
	\item[(b)]  For the particular case of probability measures, a more direct way to prove Corollary \ref{cor:Rest}  was  pointed out to us by the  referee.
It consists in observing that, as  a  consequence of Lin's theorem  (Theorem \ref{lin}), a power bounded operator $T\in L(X)$ on a Banach space $X$, is uniformly mean ergodic if and only if its restriction to a closed subspace $X_0$ (not necessarily a hyperplane) containing the $T$-invariant subspace $\overline{(I-T)(X)}$ is uniformly mean ergodic.  Indeed,  if one assumes that  $T$ is uniformly mean ergodic on $L:=\overline{(I-T)(X)}$, the ergodic decomposition of $\restr{T}{L}$ yields that   $(I-T)(X)\subseteq L=(I-T)(L)\subseteq (I-T)(X)$ so that  $(I-T)(X)$ is closed in $X$. Theorem \ref{lin} then  applies to show that $T$ is uniformly mean ergodic.
If $\mu$ is a probability measure, then $(I-\lambda_1^0(\mu))(L_1(G))\subseteq L_1^0(G)$ and this argument applies to $\lambda_1(\mu)$.
\end{itemize}	
\end{remark}

\subsection{A first characterization}
The  arguments  in the previous sections lead us to  the classic characterization of uniform mean ergodicity, summarized in Theorem \ref{firstcar:adap}. Part of that theorem, the equivalence of (2) and (4) for instance,   hold for general Markov operators, see
[21], while those involving spread-out measures are specific to convolution operators. A proof of Theorem \ref{firstcar:adap} can  be found in \cite[Section 6.3]{revuz}, see also \cite{bha73}. When the group is connected, more can be said, we refer to Theorem \ref{conexo} and the subsequent discussion.   Since, given the tools already introduced, the proof of the  theorem is rather straightforward, we have decided to   put it together using the notations of the present paper, which will be handy later on. Note that some of the references mentioned in this paragraph   deal with  convolution operators on $L_\infty(G)$ instead of $L_1(G)$. Since  convolution  operators on $L_\infty(G)$ are  adjoint to convolution operators on $L_1(G)$, both approaches yield the same classification (cf. Corollary \ref{cor:carac1:ume}, below).  This characterization will be completed in Section 5  by adding a spectral condition.

\begin{theorem}
\label{firstcar:adap}
 Let $G$ be a  compact group and $\mu\in M(G) $ be an adapted  probability measure. The following are equivalent:
 \begin{enumerate}
  \item\label{1}   $\mu$  is uniformly ergodic.
 \item\label{2} $\lambda_1(\mu)$ is uniformly mean ergodic.
 \item\label{3}   $\mu$ is spread-out.
  \item\label{4}  $\lambda_1(\mu)$ is quasi-compact.
 \item\label{5}  $\lambda_1^0(\mu)$ is quasi-compact.
 \end{enumerate}
 \end{theorem}

\begin{proof}
Proposition \ref{ume0} and  Corollary \ref{cor:Rest} prove that  (1) and (2) are equivalent.

   If $\mu$ is not spread-out, then  $\mu_{[n]}$ is a singular measure for every $n\in \N$. It follows that
  \[\norm{\lambda_1(\mu_{[n]}-\mg)}=\norm{\mu_{[n]}-\mg}\geq 1, \quad \mbox{ for every } n\in \N.\]
We deduce  from Proposition \ref{ume0} that    (2) implies  (3).

 Assume that Condition (3) holds.  There exist then $n\in \N$ and $f\in L_1(G)$ such that $\mu^n=f d\mg+\mu_s$, with $\mu_s$ singular with respect to $\mg$ and $\|\mu_s\|<1$.  This yields $\|\lambda_1(\mu^n)-\lambda_1(f\mg)\|=\|\mu^n-f\mg\|=\|\mu_s\|<1$. Since $\lambda_1(f\mg)$ is a compact operator by Theorem 4 of \cite{akemann67}, we get that $\lambda_1(\mu)$ is quasi-compact.

   If $\lambda_1(\mu)$ is quasi-compact, so will be  its restriction $\lambda_1^0(\mu)$ to $L_1^0(G)$. So  (4) implies (5).

   Finally, (5) implies (2) by Theorem \ref{yosidakaku} and Corollary \ref{cor:Rest}.
\end{proof}
%
 The implication (3)$\implies$(2) is specific to convolutions. The example in \cite[pp. 213-214]{rose81}  shows that a stationary Markov operator $P$ may satisfy (3) without being uniformly
ergodic on  $L_1$.
\begin{remark}
   By Corollary \ref{3.4} and Corollary \ref{cor:noncnorm}, Theorem \ref{firstcar:adap} actually  characterizes uniform mean ergodicity of $\lambda_1(\mu)$ for any probability measure and every locally compact group.
\end{remark}
\section{Spectral characterization of uniformly ergodic probability measures}
The deep connection between uniform ergodicity and spectral properties of operators is set forth in Theorem \ref{dunfordlin}. It is immediately clear that a uniformly mean ergodic operator cannot have 1 as an accumulation point of its spectrum. If dealing with a normal  operator  defined on a Hilbert space,  then the converse is true under mild conditions. This is the departing point of the following result proved as Theorem 4.10 in our previous paper \cite{galijorda21}.
\begin{theorem}\label{4.10}
Let $G$ be a locally compact group and let $\mu\in M(G)$ be a probability measure with $\mu^\ast \ast \mu=\mu\ast \mu^\ast$. Then $\lambda_1(\mu)$ is uniformly mean ergodic if and only if $H_\mu$ is compact and 1 is an isolated point of $\sigma(\lambda_1(\mu))=\sigma(\mu)$.
\end{theorem}

We have used here that   $\sigma(\mu)=\sigma(\lambda_1(\mu))$. This is an easy consequence of Wendel's theorem \cite{wend52} to the effect that a bounded operator $T$ on $L_1(G)$ that commutes with right translations is necessarily of the form $T=\lambda_1(\nu)$ for some $\nu\in M(G)$. In particular, when  $s\notin \sigma(\lambda_1(\mu))$, the operator $T=(\lambda_1(\mu)-sI)^{-1} $  will correspond to  the convolution operator associated to the  measure $(\mu-s \delta_e)^{-1}$.



The main objective of this section is to   remove the normality condition  $\mu^\ast \ast \mu=\mu\ast \mu^\ast$ in Theorem \ref{4.10}, and,   therefore,  to obtain a full spectral characterization of uniformly ergodic measures.
\subsection{When 1 is an accumulation of $\sigma(\mu)$, $\mu$ adapted}
Our analysis will rely on the well-known Riesz idempotent operators that spectral sets induce. We describe them below for the convenience of the reader.
\begin{theorem}[Riesz idempotent operator, Theorem 4.20 of \cite{kubr}]\label{RIO}
Let $X$ be a Banach space and $T\in\lE{X}$. Assume that $\Delta$ is a closed and open  subset   of $\sigma(T)$.
Consider a function $\psi  \colon \Lambda \to \C$  that is  analytic  on an open set $\Lambda \subseteq \C$ that contains $\sigma(T)$ and satisfies \[\psi(\Delta)=\{1\} \quad \mbox{ and } \psi\left(\sigma(T)\setminus \Delta\right)=\{0\}.\]
The function $\psi$ then defines, through Riesz Functional Calculus (see Theorem  4.7 from \cite{conway}), an  operator $E_\Delta\in \lE{X}$, with  the following properties:
\begin{alignat}{4}
   &\hspace*{-2.7cm} E_\Delta \mbox{ is a projection, }\label{RI1}\tag{1}\\
   &\hspace*{-2.7cm}\label{RI2}\mbox{If $S\in \lE{X}$ and $ST=TS$, then $E_\Delta S=SE_\Delta$,}\tag{2}\\
   &\hspace*{-2.7cm}\label{RI3} \sigma\left(T E_\Delta\right)=\Delta \cup \{0\}, \tag{3}\\
   &\hspace*{-2.7cm}\label{RI4} \sigma\left(\restr{T}{E_\Delta(X)}\right)=\Delta.\tag{4}
   \end{alignat}
\end{theorem}

We now need a  few  Lemmas.  The first one can be stated in terms of multipliers of Banach algebras. Recall that a right multiplier $M$ of a Banach algebra $\mA$ is a bounded linear operator $M\colon \mA\to \mA$ such that $M(ab)=aM(B)$.
\begin{lemma}\label{inv}
   Let $\mA$ be a Banach algebra and let  $M_1,M_2\colon \mA\to \mA$ be two commuting right multipliers of $\mA$. Assume that  $M_2$ is a projection and set $J:=M_2(\mA)$. If  $\restr{M_1}{J}$ is invertible, then there is  a right  multiplier $M^\prime \colon \mA\to \mA$ such that $\restr{M^\prime}{J}=\left(\restr{M_1}{J}\right)^{-1}$.
 \end{lemma}
 \begin{proof}
   Let $C:= \left(\restr{M_1}{J}\right)^{-1}$ and   define $\widetilde{C}\colon \mA\to \mA$ by
   \[\widetilde{C}=C\circ M_2.\] Since $\restr{\widetilde{C}}{J}=C$, it will be enough to  check that
   $\widetilde{C}$ is a right multiplier.

   Let $a, b \in \mA$. Then
   \begin{equation*}\label{eq1}
     aCM_2(b)=CM_1\left(a CM_2(b)\right)
 =C\left(a M_1CM_2(b)\right)=C\left(a M_2(b)\right).
   \end{equation*}
   It follows that
    \[ a\widetilde{C}(b)=C\left(a M_2(b)\right)=C\left(M_2(ab)\right)=\widetilde{C}(ab).\]
 \end{proof}
 The second Lemma is
 a well-known  Linear Algebra fact.
 \begin{lemma}
   \label{eigen}

   Let $T_1$ and $T_2$ be commuting bounded operators on a complex Banach space $X$. If $z\in \C$ is an eigenvalue of $T_1$ and its eigenspace $V_z=\{v\in X\: : \: T_1 v = zv\}$ is finite-dimensional, then there exists $0\neq u\in V_z$ which is an eigenvector for $T_2$.
 \end{lemma}
 \begin{proof}
    Since $T_1$ and $T_2$ commute,  $\restr{T_2}{V_z}\colon V_z\to V_z$ is   a linear operator on the  finite dimensional space $V_z$ and has therefore at least one eigenvalue.   All  eigenvectors of $\restr{T_2}{V_z}$ are  then  eigenvectors of $T_2$.
 \end{proof}%
 \begin{theorem}\label{isoso}
   Let $G$ be a compact group and $\mu\in M(G)$ be an adapted probability measure.
   The following assertions are equivalent:
   \begin{enumerate}
	\item $1$ is an isolated point of $\sigma(\mu)$.
	\item $\mu$ is spread-out.
	\item $\lambda_1(\mu)$ is uniformly mean ergodic.
   \end{enumerate}
 \end{theorem}
 \begin{proof}
 After Theorem \ref{firstcar:adap}, which shows that $(2)$ implies $(3)$, and Theorem~\ref{dunfordlin}, giving that $(3)$ implies $(1)$, it only remains to prove that $(1)$ implies $(2)$.

Since $ \{1\}$ is open and closed in $\sigma(\lambda_1(\mu))$,  there is  $0<r<1$ small enough such that $B(1,r)\cap \sigma(\mu)=\{1\}$, where $B(z,r)$ denotes the ball of radius $r$ centered in $z\in \C$. Now we consider $K:=\overline{B(0,1+\frac r2)}\setminus B(1,r)\cup \overline{B(1,\frac r2)}$. Clearly, $\sigma(\mu)\subseteq K$. Define $\psi$ on $K$ as
\[
\psi(z)=\begin{cases} 0,& \mbox{ when $z\in \overline{B}(0,1+\frac r2)\setminus B(1,r)$},  \\
1,& \mbox{ when $z\in \overline{B}(1,\frac r2)$.}\end{cases}\]  Since the function $\psi $ can be  analytically  extended to an open set $\Lambda$ containing $K$, the Riesz Functional Calculus defines an operator $E_1= \psi(\lambda_1(\mu))\in  \lE{L_1(G)}$   with the properties stated in Theorem \ref{RIO}. Properties \eqref{RI1} and \eqref{RI2} of that  Theorem  show that $E_1$ is a projection that  commutes with right translations. The classical Wendel's theorem (see \cite{wend52})  shows then that $E_1=\lambda_1(\nu)$ for some idempotent $\nu\in M(G)$.

   Let $\pi\colon G\to \mathcal{U}(\mathbb{H}_\pi)$  be an irreducible representation of $G$,  different from $\mathbf{1}$, the trivial one-dimensional representation. Assume  that $\widehat{\nu}(\pi)\neq 0$.

    If the only eigenvalue of $\widehat{\nu}(\pi)$ was 0, the characteristic polynomial of $\widehat{\nu}(\pi)$ would be $P(z)=z^n$ and the Cayley-Hamilton theorem would imply that $(\widehat{\nu}(\pi))^n= 0$. But $\nu$ is idempotent and $(\widehat{\nu}(\pi))^n=\widehat{\nu^n }(\pi)=\widehat{\nu}(\pi)$. With this and Lemma \ref{eigen}, we conclude that there  are $z_0,z\in \C$, $z_0\neq 0$, and  $\xi \in \mathbb{H}_\pi$, $\xi\neq 0$, such that  \begin{align*}
   \widehat{\nu}(\pi)\xi &=z_0 \,\xi\quad  \mbox{ and }\\
   \widehat{\mu}(\pi)\xi & =z\,\xi.
   \end{align*}
   The measure $\nu$ being idempotent we have that, actually, $z_0=1$.
      Hence, $\widehat{(\mu\ast\nu)}(\pi)\xi=z\xi$. As
 we know that $\sigma(\mu\ast \nu)=\{0,1\}$ (statement \eqref{RI3} of Theorem \ref{RIO}) and  that $1$ cannot be an eigenvalue of $\widehat{\mu}(\pi)$ (we are using here Corollary \ref{en:3} (b)), it follows from  \ref{en:4} of Theorem \ref{FS} that
    \begin{align*}
   \widehat{\nu}(\pi)\xi &= \xi \mbox{ and }\\
   \widehat{\mu}(\pi)\xi & =0.
   \end{align*}
   On the other hand,
      Statement \eqref{RI4} of Theorem \ref{RIO} shows that $0\notin \sigma \left(\restr{\lambda_1(\mu)}{\nu\ast L_1(G)}\right)$. There is therefore an operator $C\colon
   \nu\ast L_1(G)\to \nu\ast L_1(G)$ inverse to $\restr{\lambda_1(\mu)}{\nu\ast L_1(G)}$.
  We now apply Lemma \ref{inv} to
  $M_1=\lambda_1(\mu) $ and $M_2=\lambda_1(\nu)$. The conclusion of the Lemma together with Wendel's theorem provides us with a measure $\sigma\in M(G)$ such that
  \begin{equation}\label{eq2}
     \sigma\ast \mu \ast  \nu\ast f=\nu \ast f, \quad \mbox{ for every } f\in L_1(G).
  \end{equation}
    Let $f_{\pi,\xi}\in L_1(G)$ be as defined in Corollary \ref{en:3} (a).
If  we take Fourier-Stieltjes transforms in Equation \eqref{eq2} with $f=f_{\pi,\xi}$ and apply them to the representation $\pi$ and then to the vector $\xi$ (equation \eqref{eq2} is actually used in the third equality), we get:
\begin{align*}
\xi=\widehat{\nu}(\pi)\xi&=\widehat{\nu}(\pi)\widehat{f_{\pi,\xi}}(\pi)\xi\\
&=
    \widehat{ \sigma}(\pi)\widehat{ \mu}(\pi)\widehat{\nu}(\pi)\widehat{f_{\pi,\xi}}(\pi)\xi\\
    &=\widehat{ \sigma}(\pi)\widehat{ \mu}(\pi)\xi=0.
\end{align*}
This is in contradiction with our choice of $\xi$. It follows that $\widehat{\nu}(\pi)=0$ for every $\pi $ different from the trivial one-dimensional representation $\mathbf{1}$.

Since $ \widehat{\nu}(\mathbf{1})=1$ (it has to be a nonzero idempotent complex number),  we conclude that
$\widehat{\nu}(\pi)=\widehat{\mg}(\pi)$ for every irreducible representation  $\pi$ of $G$ and, hence,
(Theorem \ref{FS}, statement \ref{en:2}) that $\nu=\mg$.

We have thus that $E_1=\psi(\lambda_1(\mu))=\lambda_1(\mg)$.
%
We finally observe that, since $\mathbb{C}\setminus K$ is connected and $\psi$ can be holomorphically extended to a  an open set containing $K$, there is, by  Runge's theorem (see \cite[Corollary III.8.5]{conway}), a sequence $(\psi_n)_{n\in\N}$ of polynomials such that $\lim \psi_n=\psi$ uniformly on $K$. Having identified $E_1=\psi(\lambda_1(\mu))$ with $\lambda_1(\mg)$, we see that $\lambda_1(\mg)=\lim \psi_n(\lambda_1(\mu))$. The fact that $\lambda_1 \colon M(G) \to \lE{L_1(G)}$ is a linear isometry yields that for each $\varepsilon>0$ there is $n\in \N$ with
\[ \|\mg-\psi_n(\mu)\|<\varepsilon.\]
Now, $\psi_n$ being a polynomial, functional calculus shows that,
 for each  $n\in\N$, there are $k_n\in\N$, and  $(\alpha_{j,n})_{j=0}^{k_n}\subset\mathbb{C}$ such that $\psi_n(\mu)=\sum_{j=0}^{k_n}\alpha_{j,n}\mu^j$,
 which implies that
not all $\mu^k$ can be singular, since singular measures constitute a closed subspace of $M(G)$.
   \end{proof}

\begin{remark}
\label{R:Rajchman}
    The existence of Rajchman measures in $\mathbb{T}$ (measures whose Fourier-Stieltjes transforms vanish at infinity in $\mathbb{Z}$) that are not spread-out shows that one cannot replace ``$1$ isolated in $\sigma(\mu)$'' by ``1 isolated in $\sigma(\lambda_2(\mu))$'' (which, in the case of Abelian groups,  coincides with the \emph{natural} spectrum $\overline{\mu(\widehat{G}})$). Examples of such measures can be deduced from the proof of  Theorem 3.9 of \cite{zafran73}, see also Remark 6.4 in \cite{galijorda21}.

\end{remark}

   \subsection{Connecting the spectrum of $\lambda_1(\mu)$ with that of $\lambda_1^0(\mu)$ and $\lambda_1(\underline{\mu})$}
   \label{conn}We see in this subsection that once 1 is an accumulation point in the spectrum one of the convolution operators $\lambda_1(\mu)$, $\lambda_1^0(\mu)$ or  $\lambda_1(\underline{\mu})$, it  is  an accumulation point in the other two spectra. This leads to clean characterizations of uniform ergodicity where any of them can be used. 

 The following  notations concerning a bounded operator $T\colon X\to X$ on a Banach space $X$ will prove useful in this section.   The symbol $\sapT$ will denote  the approximate spectrum of $T$. We recall that $z\in \sapT$ whenever there exists a sequence $(x_n)$ in the unit sphere of $X$ such that $\lim \|T(x_n)-zx_n\|=0$. If $A\subset    \C$, $\partial A$ will denote the boundary of $A$  and $\acc (A)$ will stand for the set of cluster points of a set $A$.

\begin{proposition}
\label{PropAcc}
Let $T:X\longrightarrow X$ be a continuous and linear operator on a Banach space $X$ with $\|T\|=1$. Then, $1\in \acc (\sigma (T))$ if, and only if, $1\in \acc (\sapT)$.
\end{proposition}
\begin{proof}
One direction is obvious, we prove the other one. Assume that  $1\in \acc (\sigma (T))$. Then, for every $n\in \N$, there is $z_n \in \sigma (T) \bigcap \{w\in \C \: : \: d(w,1)\leq \frac1n\}\setminus\{1\}$ where $d$ here stands for the Euclidean distance in $\C$.

Let $t_n=\sup \{t\geq 0 \: : \: [z_n-ti, z_n+ti]\subseteq \sigma (T)\}$, where $[z_n-ti,z_n+t_i]$  denotes the obvious vertical segment in $\C$. Since $\sigma(T)$ is a compact set, one can choose correctly the sign so that $w_n=z_n\pm t_n i \in \partial \sigma(T)$.

Now define $x_n$ as the intersection of the line $1-\frac1n + ti$ with $\T$, the boundary of $\D$, the unit disk, so $x_n=1-\frac1n \pm \frac{\sqrt{2n-1}}{n}i$. Since $\sigma (T)\subseteq \overline\D$, we have that, by construction, $0\leq t_n\leq \frac{\sqrt{2n-1}}{n}$.
 This concludes the proof, for  $w_n\in \partial \sigma(T)\subseteq \sapT$ \cite[Chapter VII, Proposition 6.7]{conway}, and we have just seen that $\lim_n w_n = \lim_n z_n=1$.





\end{proof}

\begin{proposition}
\label{propAP}
Let $\mu\in M(G)$ be a probability measure on a non-compact group $G$. Then $\sapunocero=\sapuno$. \end{proposition}
\begin{proof}
We only have to show that $\sapuno \subseteq \sapunocero$. To that end, let $z\in \sapuno$, then there is a sequence $(f_n)_n\subseteq L_1(G)$, with $\|f_n\|_1=1$ and satisfying $\lim_n \|\mu \ast f_n-z f_n\| =0$. Since, for each $n\in \N$, $\|\lambda_1 ^0(f_n \cdot m_G)\| =1$, by Proposition~\ref{noncnorm},  there exists $g_n\in L_1^0(G) $, with $\|g_n\|_1=1$ such that $\frac12 < \|f_n \ast  g_n\|_1\leq 1$.

We conclude by observing that $ (f_n\ast g_n)/\norm{f_n\ast g_n}_1\in L_1^0(G)$ and that
\begin{align*}
 \lim_{n\to \infty} \left\|\mu \ast\frac{f_n\ast g_n}{\|f_n\ast g_n\|_1}-z\cdot\frac{f_n*g_n}{\|f_n*g_n\|_1}\right\|_1&\leq 2\lim_{n\to \infty} \|\mu *f_n*g_n-z\cdot f_n*g_n\|_1
 \\
 &\leq 2\lim_{n\to \infty} \|\mu *f_n-z\cdot f_n\|_1=0.
\end{align*}
\end{proof}

%
%
%
%
%
%
\begin{corollary}\label{cor:propAP}
Let $G$ be a locally compact group and let $\mu\in  M(G) $ be a probability measure.
\begin{itemize}
\item[(a)] 1 is isolated in $\sigma(\lambda_1(\mu))$ if and only if $1$ is not an accumulation point of $\sigma(\lambda_1^0(\mu))$.
\item[(b)]  1 is isolated in $\sigma(\lambda_1(\mu))$ if and only if 1 is isolated in $\sigma(\lambda_1(\underline{\mu}))$.
\end{itemize}
\end{corollary}
\begin{proof}
We  first prove (a).   Our argument will be different depending on whether $G$ is  compact or not.

  If $G$ is compact, the map $f\mapsto (f-\int fdm_G) \oplus \int fdm_G$ establishes  a linear isometry between $L_1(G)$ and $L_1^0(G)\oplus \ \C$. Moreover, since  both the closed hyperplane $L_1^0(G)$  and its complement are  $\lambda_1(\mu)$-invariant, we can decompose $\lambda_1(\mu)$ as $\lambda_1(\mu)=\lambda_1^0(\mu)\oplus \text{Id}$. It follows that $\sigma(\lambda_1^0)(\mu) \cup\{1\}=\sigma(\lambda_1(\mu))$, which certainly implies the result.

  If   $G$ is not compact,  Propositions \ref{PropAcc} and  \ref{propAP} together yield directly the result.

We now prove Statement (b).
Since $M(H_\mu)\subseteq M(G)$ as unitary Banach algebras, it follows
$$\sigma(\lambda_1(\mu)) =\sigma(\mu)\subseteq \sigma(\underline{\mu})=\sigma(\lambda_1(\underline{\mu})).$$

\noindent Again, we just  have to show that if 1 is an accumulation point of $\sigma(\lambda_1(\underline{\mu}))$ then 1 also accumulates in $\sigma(\lambda_1(\mu))$. We can deduce from Proposition \ref{PropAcc} that it suffices to show that $\sap(\lambda_1(\underline{\mu}))\subseteq \sapuno$. Let  $\lambda\in \sap(\lambda_1(\underline{\mu}))$. There is $f_n\in L_1(H_\mu)$ such that $\|f_n\|=1$ and $\lim_n \|\mu*f_n-\lambda f_n\|=0$.  Now,  we observe that, for every $n$, $f_n\dmh\in  M(G)$ and $1=\norm{f_n}=\norm{\lambda_1(f_n\dmh)}$, there will be therefore
 $g_n\in L_1(G)$ , $\|g_n\|=1$, such that $\|\lambda_1(f_n)(g_n)\|=\|f_n\ast g_n\|>\frac12$. Choosing $h_n:=(1/\|f_n\ast g_n\|)f_n\ast g_n\in L_1(G)$ we get
 $$\lim_n \|\mu*h_n-\lambda h_n\|=0,$$

 \noindent and $\lambda\in \sapuno$.


\end{proof}
\subsection{The spectral characterizations}
The results  of the previous subsections can now be put to work.


\begin{theorem}
\label{car:adap}
 Let $G$ be a  compact group and $\mu\in M(G) $ be an adapted  probability measure. The following are equivalent:
 \begin{enumerate}
 \item   $\mu$  is uniformly ergodic.
  \item $\lambda_1(\mu)$ is uniformly mean ergodic.
  \item   $\mu$ is spread-out.
  \item  $\lambda_1(\mu)$ is quasi-compact.
 \item  $\lambda_1^0(\mu)$ is quasi-compact.
 \item  $1$ is an isolated point of $\sigma(\mu)$
\item  $1$ is not  an accumulation point of $\sigma(\lambda_1^0(\mu))$.
 \end{enumerate}
 \end{theorem}
 \begin{proof}
   Statements (1)-(5) have already been shown to be equivalent, Theorem \ref{firstcar:adap}. Statement (2) implies Statement (6)  by Theorem \ref{dunfordlin} together with the equality $\sigma(\lambda_1(\mu))=\sigma(\mu)$. Statement (6) implies Statement (3), by Theorem \ref{isoso}. Statements (6) and (7) are equivalent by Corollary \ref{cor:propAP} (a).
 \end{proof}

 We can also characterize uniform mean ergodicity of $\lambda_1(\mu)$ for any $\mu \in M(G)$, in terms of $\mu$ itself.
 Recall that  we denote by $\underline{\mu}$, the measure $\mu$ seen as a measure in $ M(H_\mu)$.
\begin{theorem}\label{carac1:ume}
  Let $G$ be a locally compact group and let $\mu \in M(G)$ be a probability measure.  The following statements are equivalent:
  \begin{enumerate}
        \item  $\lambda_1(\mu)$ is a uniformly mean ergodic operator.
        \item  $\lambda_1^0(\mu)$ is a uniformly mean ergodic operator.
         \item  The measure  $\underline{\mu}$ is uniformly ergodic.

    \item $H_\mu$ is compact and 1 is an isolated point of $\sigma(\mu)$.
      \end{enumerate}
\end{theorem}
\begin{proof}
The equivalence of the first three statements follows from Theorem \ref{car:adap} and Corollary \ref{cor:Rest}. If Corollary \ref{cor:noncnorm} is taken into account,  Theorem \ref{car:adap} also shows that these  three statements are equivalent to $H_\mu$ being compact and $1$ being isolated in $\sigma(\lambda_1(\underline{\mu}))$. We conclude with  Corollary \ref{cor:propAP} (b).
\end{proof}

 As already remarked after Theorem\ref{rose81}, the measure $\mu^\ast$ needs not be ergodic when $\mu$ is. Furthermore, M. Rosenblatt  \cite{rose71} produces an example of a stationary (irreversible) Markov chain $P$ that  defines an operator $T$ which  is uniformly mean ergodic on $L_1$ while $T$ is not uniformly mean ergodic on $L_\infty$, so $T^\ast$
is not uniformly mean ergodic on $L_1$.
 Corollary \ref{cor:carac1:ume} shows that adjoints behave better  in our setting.



\begin{corollary}\label{cor:carac1:ume}
  Let $G$ be a locally compact group and $\mu \in M(G)$ be a probability measure.   The following statements are equivalent to the Statements (1)-(4) of Theorem \ref{carac1:ume}.
  \begin{enumerate}
  \addtocounter{enumi}{4}
        \item $\lambda_\infty(\mu)$ is mean ergodic.
        \item $\lambda_\infty(\mu^\ast)$ is mean ergodic.
        \item $\lambda_1(\mu^\ast)$ is uniformly mean ergodic.
    \end{enumerate}
\end{corollary}
\begin{proof}

We  need to recall for this proof that
the adjoint  operator of $\lambda_1(\mu)$ is the operator $\lambda_\infty(\mu^\ast):L_\infty(G)\longrightarrow L_\infty (G)$ \cite[Theorem  20.23]{hewiross1}.
Since inverse sets of Haar-null sets are Haar-null,  and $(\mu^\ast)^n=(\mu^n)^\ast$,  \cite[Theorem 20.22]{hewiross1}), $\underline{\mu}$ will be spread-out if and only if
$\underline{\mu}^\ast$ is. Theorem \ref{carac1:ume} then shows that $\lambda_1(\mu)$ is uniformly mean ergodic if and only if $\lambda_1(\mu^\ast)$ is. Since the adjoint of a uniformly mean ergodic operator is  uniformly mean ergodic as well and $\lambda_\infty(\mu)$ is the adjoint of  $\lambda_1(\mu^\ast)$, the equivalence of Statement (5) here and Statement (1) of Theorem \ref{carac1:ume} follows. To conclude the proof we only need Lotz's theorem \cite[Theorem 5]{lotz}, which ensures that a power bounded operator on $L_\infty(G)$ is uniformly mean ergodic whenever it is ergodic.
\end{proof}

We also note the following \emph{abelian} consequence of Theorem \ref{isoso}.

\begin{corollary}
Let $G$ be a locally compact  \emph{abelian} group. The following are equivalent.
\begin{enumerate}
  \item  $\lambda_1(\mu)$ is uniformly mean ergodic.
  \item 1 is isolated in $\sigma(\mu)$.
 \end{enumerate}
\end{corollary}
\begin{proof}
Theorem \ref{carac1:ume} shows that (1) implies (2). By the same characterization, the  converse statement will be proved if we show that 1  isolated in $\sigma(\mu)$ implies that $H_\mu$ is compact.

If $H_\mu$ is  not compact, then $\widehat{H_\mu}$ is  not discrete and one can find a net $\{\chi_\alpha\}_\alpha\subseteq\widehat{H_\mu}\setminus\{\mathbf{1}\} $ that converges to the trivial character $ \mathbf{1}$.  The net $\{\widehat{\underline{\mu}}(\chi_\alpha)\}_\alpha$ is then contained in $\T\setminus\{1\}$, by Corollary \ref{en:3},  and converges to 1. It follows that 1 is not isolated in $\sigma(\lambda_1(\underline{\mu}))$. Corollary \ref{cor:propAP}.(b) then shows that 1 is not isolated in $\sigma(\mu)$.
\end{proof}

\begin{remark}
As we have observed, uniform ergodicity of $\lambda_1(\mu)$ is a strong property. Indeed, uniformly ergodic stationary Markov operators on $L_1$ are  uniformly mean ergodic on each $L_p$, $1<p<\infty$, \cite[Theorem 1, p. 211]{rose71}. The converse does  not hold, the adjoint $T^\ast$ of the operator $T$ constructed in  \cite[pp. 213-214]{rose81}  is uniformly mean ergodic in $L_2$, since $T$ is, but is not uniformly
mean ergodic in $L_1$. Careful constructions involving  Rachjman measures provide examples of convolution operators that are  uniformly mean ergodic on  $L_2$ but are not uniformly ergodic in  $L_1$. Constructions of that sort can be found  in Theorem 5.6 of \cite{brownkarawill82} (where Rajchman measures whose spectrum is the whole disk are  constructed) and in Theorem 6.2 of \cite{cohlin22} (where Rajchman measures that  are not spread-out are  constructed).

It may be worth to remark that, although the  operator  $\lambda_2(\mu) $ can be uniformly mean ergodic without $\lambda_1(\mu)$ being so,  the proof of the preceding Corollary shows that, still, for $\lambda_2(\mu)$ to be uniformly mean ergodic it is necessary that $H_\mu$ is compact. \end{remark}

\begin{remark}
The authors of \cite{mustatopal21} observe  that Host and Parreau's characterization  of the measures $\mu\in M(G)$ such that $\mu\ast L_1(G)$ is closed in $M(G)$ \cite{hostparr78} implies that   $\lambda_1(\mu)$ is uniformly mean ergodic if and only if $\mu$ is of the form $\mu=\delta_e+\lambda\ast \theta$, with $\lambda$ an invertible measure and $\theta$ an idempotent measure.  Hence, Theorem \ref{car:adap} and Theorem \ref{carac1:ume} together show that this description for probability measures is equivalent to $H_\mu$ being compact and $\underline{\mu}$ being spread-out.
\end{remark}
\subsection{Uniform convergence of  powers of convolution operators}

Our results on uniform ergodicity  can also be translated to uniform complete mixing and, more generally, to the convergence of the powers of $\lambda_1(\mu)$ and $\lambda_1^0(\mu)$.

We first characterize uniformly completely mixing measures. In this characterization, the equivalence between (1) and (3) is Theorem 4.1 of \cite{anougatz04}. Since condition (2) is trivially equivalent to (1) under the considered hypothesis,   only the equivalence between (4) and (2) needs to be proved. For the sake of completeness, we  also offer here a short alternative proof  of the equivalence between the first three statements based on Theorem \ref{car:adap} and classic ergodic arguments.

\begin{theorem}\label{portmantCM}
Let $G$ be  a compact group and let $\mu \in M(G)$ be an adapted probability measure. The following assertions are equivalent.
\begin{enumerate}
        \item  The sequence $(\lambda_1(\mu^n))$ is norm convergent.
          \item  The measure  $\mu$ is uniformly completely mixing.
          \item   $\mu $ is spread-out  and  strictly aperiodic.

\item  $1$ is isolated in $\sigma(\mu)$ and $\sigma(\mu)\cap \T=\{1\}$.
\item $\lambda_1(\mu)$ is uniformly mean ergodic and $\|\lambda_1(\mu^{n+1})-\lambda_1(\mu^n)\|\to 0$.
    \end{enumerate}
\end{theorem}

\begin{proof}

 Since $G$ is compact, $\lambda_1(\mu)=\lambda_1^0(\mu)\oplus \mathrm{Id}$, hence Proposition \ref{ume0} shows that (1) and (2) are equivalent.
 Condition (2) implies that $\lambda_1^0(\mu)$ is uniformly mean ergodic and $\mu$ is completely mixing. Uniform ergodicity   implies that $\mu$ is spread-out (Theorem  \ref{car:adap}). Completely mixing and Lin-Wittmann's Theorem \ref{adp-spbasico}.(b) yield that $\mu$ is strictly aperiodic. Hence (2) implies (3).

 If we assume (3), then $\lambda_1^0(\mu)$ is quasi-compact by Theorem \ref{car:adap}, and Kawada-Ito's theorem \ref{adp-spbasico}.(3) yields that $(\lambda_1^0(\mu^n))$ is SOT convergent to 0,   hence $\sigma_p(\lambda_1^0(\mu))\cap \T=\emptyset$. Yosida-Kakutani's theorem \ref{yosidakaku} implies then that $(\|\lambda_1^0(\mu^n)\|)$ is convergent to 0, hence we have (2). Statements (1)-(3) are thus shown to be equivalent.

   Assume that  the first three equivalent conditions hold. In a Banach algebra the elements whose sequence of powers is convergent to 0 are those whose spectral radius is smaller than 1, hence (2)  implies $r(\lambda_1^0(\mu))<1$. Since $\sigma(\mu)=\sigma(\lambda_1(\mu))=\sigma(\lambda_1^0(\mu))\cup \{1\}$, (4) necessarily holds.

    Suppose conversely that  (4) holds. Then $\mu$ is uniformly ergodic (Theorem \ref{car:adap}) and    $1\notin \sigma(\lambda_1^0(\mu))$, for otherwise  also $1\in \sigma(\lambda_1^0(\mu_{[n]}))$ for each $n\in\N$, against uniform ergodicity. Using again that $\sigma(\mu)=\sigma(\lambda_1(\mu))=\sigma(\lambda_1^0(\mu))\cup \{1\}$ we conclude that $r(\lambda_1^0(\mu))<1$, which is equivalent to (2).


    Condition (5) being equivalent to condition (4) is the theorem of Katznelson and Tzafriri \cite{katztza86}.
\end{proof}

\begin{remark} The previous  theorem, together with Corollary \ref{noncnorm},  provides a  solution to the uniform version of the completely mixing: a uniformly ergodic and strictly aperiodic
measure is necessarily uniformly completely mixing.
\end{remark}

\smallskip
If we do not assume $\mu$ to be adapted, Theorem \ref{portmantCM} can be stated as follows.
\begin{theorem}
\label{iterates}
 Let $G$ be a locally compact group and let $\mu \in M(G)$ be a probability measure. The following statements are equivalent:
  \begin{enumerate}
        \item  $(\lambda_1(\mu^n))$ is convergent in $(\lE{L_1(G)},\|\cdot\|)$.
        \item  $(\lambda_1^0(\mu^n))$ is convergent in $(\lE{L_1^0(G)},\|\cdot\|)$.
        \item $\lambda_1^0(\mu)$ is uniformly mean ergodic and $(\lambda_1^0(\mu^n))$ is convergent in $\lE{L_1^0(G)}$ endowed with the strong operator topology.
        \item $H_\mu$ is compact, $\underline{\mu}$ is strictly aperiodic  and 1 is isolated in $\sigma(\mu)$.
      \end{enumerate}

\begin{proof}
(1) implies (2) and (2) implies (3) are trivial. Assertion (3) implies that $\underline{\mu}$ is uniformly ergodic by Theorem \ref{carac1:ume}. This implies that  $\sigma_p(\lambda_1^0(\underline{\mu}))\cap \T\subseteq \{1\}$ and, by  Theorem \ref{car:adap}, that  $\lambda_1^0(\underline{\mu})$ must be  quasi-compact. Theorem \ref{yosidakaku} proves then that $\underline{\mu}$ is uniformly completely mixing. Statement (4) then follows from Theorem \ref{portmantCM}, via Corollary \ref{cor:noncnorm},  and Corollary \ref{cor:propAP}. This latter Corollary shows that we can get   Statement (1) from Statement (4) by applying  Theorem \ref{portmantCM} to $\underline{\mu}$  and then Corollary \ref{3.4}.
\end{proof}

\end{theorem}
 If  $G$ is connected,  the condition of strict aperiodicity can be dropped from Theorem \ref{portmantCM}.

\begin{theorem}
\label{conexo}
Let $G$ be a  locally compact  \emph{connected} group and $\mu \in M(G)$, be a probability measure. Then the following statements are equivalent:
\begin{enumerate}
    \item $\mu$ is uniformly completely mixing.
    \item $\mu$ is uniformly ergodic.
    \item $G$ is compact and $\mu$ is spread-out.
    \end{enumerate}
\end{theorem}
\begin{proof}
It is  trivial that (1) implies (2). Conditions    (2) and (3) are equivalent by Corollary \ref{cor:noncnorm} and Theorem \ref{firstcar:adap}. We see that the conjunction of (3) and (2) implies (1). Let $\mu$ be a uniformly ergodic  measure a compact group $G$. By Theorem \ref{portmantCM} we only need to show that $\mu$ is strictly aperiodic. We proceed by contradiction and  assume that there are a proper normal closed subgroup $H$ and a point $x\in G$ such that $S_\mu \subseteq xH$. Note that $\mg (H) =0$, for else it would be open and, bearing in mind that $G$ is connected, we would have $H=G$. Now we have
\[
\bigcup_{1\leq j,k\leq n} S_\mu ^{-j} S_\mu ^k \subseteq \bigcup _{j=-n} ^n x^j H,
\]

\noindent for each $n\in \N$. Since $\mg(H)=0$ we can apply Proposition~\ref{caractporconjuntos}, to  conclude that $\mu$ is not uniformly ergodic.
\end{proof}

\begin{remark}
 If we add $G$ compact to the hypothesis of the previous Theorem  (a condition that, as it turns out, is necessary), the resulting statement  is completely equivalent  to  Theorem 3 of \cite{bha73} where  Bhattacharya proves that the sequence $(\mu^n)$ converges  to $\mg$ for every non-singular probability measure $\mu$  on a connected compact group $G$.  This statement is formally weaker than Theorem \ref{conexo}.
If one only assumes that $\mu$ is spread out as is done in Theorem \ref{conexo},
\cite[Theorem 3]{bha73} would just show that $\lim_n \norm{\mu^{k n}-\mg}=0$ for some $k\in \N$. But, as a matter of fact, from here one readily deduces that $\mu$ is uniformly completely mixing. It suffices to factorize,  for any $n\in \N$, $\mu^n=\mu^{km_n}\ast \mu^{j_n}$, for some $0\leq j_n< n$ and $m_n\in \N$ and recall that $\mg\ast\nu=\mg$ for any probability measure $\nu$.

Other  proofs of this  can be found in literature, see for instance  \cite[Chapter 6, Exercise 3.19]{revuz} or  \cite[Corollary 4.2, Remark (5)]{anougatz04}.

\end{remark}

\begin{example}
  Let $A$ be an arc 0 in $\T$ of length less than 1/2, centered at 1. Then $\mu=\frac{1}{\m{\T}(A)}\Cf_{A}\m{\T}$ is uniformly completely mixing, even if $\norm{\lambda_1^0(\mu) }=1$, by Lemma \ref{prop:apro}.
  \end{example}

\begin{remark} The above  example shows that for an adapted probability measure $\mu$ on a compact group , contrary to what happens with $\lambda_1(\mu)$,  the inequality $r(\lambda_1^0(\mu))<\|\lambda_1^0(\mu)\|$ is possible. The inequality $r(\lambda_1^0(\mu))<1$ for a probability $\mu\in M(G)$ is actually equivalent, for any locally compact group $G$,  to $\mu$ being uniformly completely mixing by the very definition and basic theory of Banach algebras. Hence, the inequality $r(\lambda_1^0(\mu))=\|\lambda_1^0(\mu)\|=1$ is fulfilled by any
  probability $\mu\in M(G)$ which is not uniformly completely mixing.
\end{remark}

If $G$ is abelian but not connected, measures failing to satisfy Theorem~\ref{conexo} can always be constructed. We thus have the following characterization of connected  groups in the category of compact abelian groups.

\begin{theorem}
\label{ab}
Let $G$ be a compact abelian group. The following assertions are equivalent:
\begin{enumerate}
\item G is  connected.
\item Every uniformly ergodic probability measure on $G$  is  uniformly completely mixing.
\end{enumerate}
\end{theorem}
\begin{proof}
  After Theorem \ref{conexo} we only have to show that Statement (2) implies Statement (1).  Suppose to that  end that $G$ is not connected. Then $G$ has a maximal proper normal subgroup $H$ which is open (since $\widehat G$ is not torsion-free \cite[24.25]{hewiross1}, there must exist $\chi\in \widehat G$ such that $\langle \chi\rangle$ has a prime number of elements, then  $H=\langle \chi\rangle^\perp$).

 Let then  $\mu= \Cf_{xH}\mg$, for some $x\notin H$. It is easy to see that $H\subsetneq \langle xH \rangle$, while the maximality of $H$ shows that $\overline{H_\mu}=\overline{\langle xH \rangle}=G$.
  The  measure $\mu$ is hence adapted and spread-out, but it is not strictly aperiodic by construction. Theorem \ref{firstcar:adap} and Theorem \ref{adp-spbasico} then suffice to show that $\mu$ is  uniformly ergodic probability  but it   is not uniformly completely mixing.
\end{proof}
 Commutativity is essential in the above example. If $G$ is a simple group or, more generally, if $G$ admits no nontrivial continuous  characters (see \cite[Lemma 3.4]{linwitt94}), then every adapted measure is strictly aperiodic.  The theorem is no longer true either, if the condition \emph{uniform} is dropped. If $x\in \T$ is not a root of unity, then $\delta_x\in M(\T)$ is ergodic but not strictly aperiodic.
%


%

\subsection{A characterization of quasi-compactness of convolution operators}


Next theorem should be compared to Theorem 2 from \cite{Lin75} and Theorem 3.4 of \cite{revuz}, where the equivalence between (1) and (2) below is proved for a more general class of operators.

\begin{theorem}\label{carac:quasi}
Let $G$ be a locally compact group and let $\mu$ be a probability measure. Let $T$ stand for $\lambda_1(\mu)$, $\lambda_\infty(\mu)$ or $\lambda_1^0(\mu)$. The following assertions are equivalent:
\begin{enumerate}	
\item	The operator $T$ is quasi-compact.
\item $(T_{[n]})_n$ is norm convergent to a finite dimensional projection.
\item $G$ is compact, $H_\mu$ is open in $G$ and $\underline{\mu}$ is spread-out.
\end{enumerate}
\end{theorem}

\begin{proof}
That Statement (1) implies Statement (2) is the Yosida-Kakutani theorem, Theorem \ref{yosidakaku}.

Let us see that (2) implies (3). If (2) holds, then $H_\mu$ is compact, by Corollary \ref{cor:carac1:ume}. We now  proceed by contradiction and  suppose that $G$ is either  noncompact or that  $G$ is compact and $H_\mu$ is not open in $G$. In both cases, we can choose a sequence $(x_k)\subseteq G\setminus H_\mu$ such that $(H_\mu x_k)$ is a  disjoint sequence of  compact subsets of $G$. We fix $n\in\N$. One can get a compact  neighbourhood of the identity $U$  such that $H_\mu x_iU\cap H_\mu x_jU=\emptyset$ whenever $1\leq i<j\leq n$. To find $U$, first we choose inductively compact neighbourhoods of the identity $U_i$, $1\leq i\leq n$ such that $H_\mu x_i\cap H_\mu x_lU_l=\emptyset$ when $1\leq l<i\leq n$ and $H_\mu x_iU_i\cap  H_\mu x_j=\emptyset$ for $1\leq i<j\leq n$. The set $U:=\cap_{i=1}^{n}U_i$ satisfies the required condition. We observe  that if $f_i$ is the characteristic function of $x_iU$ then $\supp (\mh\ast f_i)=H_\mu x_i U$ for each $i\neq j\leq n$. Hence the dimension of $\lambda_1(\mh)(L_1(G))$ and the dimension of $\lambda_\infty(\mh)(L_\infty(G))$ are at least $n$. By observing  how the convolution operator acts on $f_{ij}=f_i-f_j$, $i\neq j$, we get also that the dimension of $\lambda_1^0(\mh)(L_1^0(G))$ is at least $n$. Since $n$ is arbitrary, we get that, in all the considered cases, the projection limit of $T_{[n]}$ is not compact.

We see finally that  Statement (3) implies Statement (1). If $G$ is compact and $\underline{\mu}$ is spread-out,  then there exists $n$ such that $\mu^n$ and  $\mathrm{m}_{H_\mu}$ are not singular. If we assume in addition $H_\mu$ to be open in $G$, then  $\mu^n$ and $\mg$ are not singular either.
Hence, $\lambda_1(\mu)$ is quasi-compact.

 \end{proof}

\begin{example}
\label{fin}
Let $G=\T\oplus K$, $K$ compact. Let $f\in L_1(\T)$, $f\geq 0$ and $\int_\T f dm_\T =1$, and consider $\mu:=f\m{{\T}}\oplus \delta_e$. The operator $\lambda_1^0(\mu)$ is then quasi-compact in case $K$ is finite, but it is not if $K$ is infinite. In both cases, $\lambda_1^0(\mu)$ is uniformly mean ergodic, even $(\lambda_1^0(\mu^n))$ is  convergent in the norm topology.
\end{example}

\begin{remark}
\leavevmode
	\begin{itemize}
		\item[(a)]	The equivalence between conditions (2) and (6) (or
(5) and (7)) in  Theorem \ref{car:adap} and the equivalence between conditions (4) and (5) in Theorem \ref{portmantCM} do not hold for general operators. The following example was kindly provided by one of the referees. Let $V$ be the Volterra operator on $L_2([0,1])$ and $T=I-V$. It is known that  $\sigma(V)=\{0\}$, and hence that $\sigma(T)=\{1\},$  (and then $\|T^{n+1}-T^n\|\to 0$ by Katznelson and Tzafriri \cite{katztza86}), and that $T$ is power bounded,
 see e.g. \cite{allan97}. It follows that $T^n\to 0$ in the strong operator topology. But, given that      $1\in \sigma(T_{[n]})$ for each $n\in\N$,  $\|T_{[n]}\|$ does not converge to $ 0$ (and hence  $\|T^n\|$ does not converge to 0, either).

\item[(b)]In Theorem \ref{iterates},  we have seen that for $T=\lambda_1^0(\mu)$, the sequence $(T^n)$ is norm convergent  if and only if $T$ is uniformly mean ergodic and the sequence $(T^n)$ is SOT convergent. Once again, one could wonder if this  statement holds for any operator. As in the previous remark, it does not.  
Let $a=(a_n)$ be a sequence of numbers such that $\lim a_n=-1$  and $|a_n|<1$ for all $n\in\N$. The multiplication operator $M_a:l_1(\N)\to l_1(\N)$, $(b_n)\mapsto (a_nb_n)$ satisfies $\|M_a\|=1$ and $1\notin\sigma(M_a)=\{a_n:\ n\in\N\}\bigcup \{-1\}$, hence $M_a$ is uniformly mean ergodic by Theorem \ref{lin}. Moreover, $M_a$ is easily seen to be SOT-convergent to 0. However $\|M^n_{a}\|=1$ for all $n\in \N$.

 When $T$ is quasi-compact, norm convergence of $T^n $ does follow from uniform mean ergodicity and SOT convergence of $T^n$. This is because     $\sigma_p(T)\cap \T\subseteq \{1\}$ whenever $(T^n)$ is a SOT convergent sequence and  Theorem \ref{yosidakaku} applies. In Theorem \ref{carac:quasi} and Example \ref{fin} we see that, in general, uniform mean ergodicity and SOT convergence of the sequence $(\lambda_1^0(\mu^n))$ does not imply quasi-compactness of $\lambda_1^0(\mu)$.
\end{itemize}
\end{remark}

\subsection*{Acknowledgements}
The authors would like  to express their gratitude to the  anonymous referees, who read the paper carefully, and provided us with comments, references and suggestions which certainly have contributed to improve this work.

The  contribution of J. Galindo to this article is part of the grant PID2019-106529GB-I00 funded by
 MCIN/AEI/10.13039/501100011033.

The research of E. Jordá was partially supported by the project PID2020-119457GB-100 funded by MCIN/AEI/10.13039/501100011033 and by “ERDF A way of making Europe”. The research of E. Jordá is also partially supported by GVA-AICO/2021/170.

A. Rodr\'iguez-Arenas acknowledges the support of Acci\'o 3.2 POSDOC/2020/14 of UJI. The research was supported partially by the grant PID2019-106529GB-I00 funded MCIN/ AEI /10.13039/501100011033.

%
\end{document}